\documentclass[11pt]{amsart}

\usepackage{graphicx}
\usepackage{amssymb}
\usepackage{amscd}

\newtheorem{thm}{Theorem}[section]
\newtheorem{prop}[thm]{Proposition}
\newtheorem{lem}[thm]{Lemma}

\newtheorem*{thm*}{Theorem}
\theoremstyle{definition}

\newtheorem*{defn}{Definition}

\newtheorem*{eg}{Example}

\newcommand{\bq}{\mathbb{Q}}
\newcommand{\bz}{\mathbb{Z}}
\newcommand{\br}{\mathbb{R}}
\newcommand{\bc}{\mathbb{C}}

\newcommand{\fred}{S^2\times S^2}
\newcommand{\lan}{\left\langle}
\newcommand{\ran}{\right\rangle}

\newcommand{\homeo}{\cong}

\newcommand{\fix}{\operatorname{Fix}}

\newcommand{\arf}{\operatorname{Arf}}

\title{Concordance of $\bz_p\times\bz_p$ actions on $S^4$ }
\author{Michael McCooey}\date{\today}
\begin{document}
\maketitle

\section{Introduction}

The main result of ~\cite{MM3} was that if $M$ is a simply-connected four-manifold admitting an effective,  homologically trivial, locally  linear action by $G=\bz_p\times\bz_p$, where $p$ is prime, then\footnote{with one exception: pseudofree actions on $\bz_3\times\bz_3$ on homotopy $\bc P^2$s.}
$M$ is equivariantly homeomorphic to a connected
sum of standard actions on copies of $\pm\bc P^2$ and $\fred$ with a possibly non-standard action on $S^4$. 
In this note we further examine these non-standard actions on the sphere. We describe some constructions arising from counterexamples to the generalized Smith Conjecture and then consider the classification of actions up to concordance. An analysis of singular sets and quotient spaces, combined with an application of the results of~\cite{CS2}, allows us to  prove:

\begin{thm} \label{maintheorem} A locally linear action of $\bz_p\times\bz_p$ on $S^4$ is topologically concordant to a linear action if and only if 
a certain Kervaire-Arf invariant $c\in \bz_2$ vanishes. Thus after normalizing rotation numbers, there are at most two concordance classes of $\bz_p\times\bz_p$ actions on $S^4$.
\end{thm}

The invariant $c$ is well-behaved with respect to two connected-sum constructions.
We finish with a discussion of the existence of actions on which $c$ is nontrivial, which is still an open question.


\section{Warmup: $\bz_p$ actions on $S^4$}

Smith Theory~\cite{SmithTransformations1, Bredon} shows that when a cyclic group of prime order $p$ acts on a sphere, the fixed-point set is again a $\bz_p$-homology sphere. 

The classical Smith Conjecture (whose solution  is now a theorem which is old enough to drink beer:~\cite{SmithConj}) goes further, stating the fixed-point set of a tame cyclic group action on $S^3$ is either empty, or an \emph{unknotted} $S^1$.  But the natural generalization of the Smith conjecture to higher dimensions is false: Giffen~\cite{Giffen}, Gordon~\cite{Gordon2}, and Sumners~\cite{Sumners} all construct counterexamples in spheres of dimension 4 and higher. Along with ~\cite{MM3}, these counterexamples motivate our study.

Suppose $\langle g\rangle \cong\bz_p$ acts on $S^4$, locally linearly and preserving orientation. $\fix(g)$ is then a $\bz_p$-homology sphere of dimension $0$ or $2$; since $\fix(g)$ is also a submanifold,  it must be either $S^0$ or $S^2$. Both are possible, but let us assume $\fix(g)\cong S^2$.
Let $N(\fix(g))$ be an equivariant regular neighborhood (see ~\cite[9.3A]{FQ}).
Then $S^4\setminus N(\fix(g))$ has the same boundary, and by Alexander duality, the same integral homology, as $S^1\times B^3$.

\begin{lem}
$S^4/\langle g\rangle$ is homeomorphic to $S^4$.
\end{lem}

\begin{proof}

First note that near a fixed point, the quotient map takes the form $(z,w)\mapsto (z,w^p)$ in complex coordinates. 
Thus the quotient orbifold is topologically a manifold.

Let $X=S^4/\langle g\rangle$. Then $\pi_1(S^4-N(\fix(g)))$ is a normal subgroup of $\pi_1(X\setminus(N(\fix(g))/\lan g \ran))$ of index $p$.
Filling in $N(\fix(g))/\langle g\rangle$ kills a meridian $\mu$ of $\fix(g)$. This in particular kills $\mu^p$. But $\mu^p$
normally generates $\pi_1(S^4-N(\fix(g)))$, and each $\mu^i$ represents one of its $p$ cosets. It follows that
$X$ is simply connected. 

Finally, the transfer map defines an isomorphism between $H_2(X ;\bq)$ and $H_2(S^4;\bq)^{\bz_p}=0$, so $X$ is a rational, and hence integral, homology sphere. It follows
from the four-dimensional topological Poincar\'e conjecture that $X\cong S^4$.

\end{proof}

\begin{defn}
We will say actions $\psi_0$ and $\psi_1$ of a group $G$ on a manifold $M$ are \emph{concordant} if
there is a locally linear action $\Psi$ of $G$ on $M\times I$ such that $\Psi|_{M\times\{i\}}=\psi_i$, and such that 
for each $g\in G$, $\fix(g, S^4\times I) \cong\fix(g, S^4)\times I$. 

\end{defn}

\noindent{\bf Remarks:}

\begin{enumerate}
\item{This stratified notion of concordance is \emph{a priori} stronger  than an unstratified version would be. For example, Bredon~\cite[1.7]{Bredon} describes involutions of $S^5$ whose fixed-point sets are lens spaces, and puncturing such an example twice along its fixed set yields a concordance from a standard $\bz_2$ action on $S^4$ to itself whose fixed stratum is not a cylinder.}

\item{Rotation numbers are carried along tubular neighborhoods of the singular strata, so it follows immediately that concordant actions must have the same rotation numbers.}
\end{enumerate}


We begin with an observation which is essentially an equivariant version of a theorem of Kervaire, and sketch a proof. Details of the argument are fleshed out in its generalization to $\bz_p\times\bz_p$.

\begin{thm}\label{rankonecase}
Every $\bz_p$ action on $S^4$ which fixes a $2$-sphere is topologically concordant to a linear action.
\end{thm}

\begin{proof}
 (Sketch) Consider the quotient space pair  $(S^4, \fix(g))/\langle g\rangle$. According to Kervaire ~\cite{Kervaire1965}, 
all even-dimensional knots are slice, so the pair  $(S^4, \fix(g))/\langle g\rangle$ bounds a ball pair $(B^5, B^3)$.


Cappell and Shaneson's
homology surgery groups~\cite{CS2} measure the obstruction to surgering $(B^5 -N(B^3))$ rel boundary to make it $\bz[\bz_p]$-homology equivalent to a standard pair  $(B^5-N(B^3_{\text{std}}))$. The obstruction groups vanish in this case, so the $\bz_p$-cover $\widetilde{B^5-N(B^3)}$ may be assumed to be a homology $S^1\times B^4$ with fundamental group $\bz\times\bz$. Filling in $S^2\times I$ makes the covering simply connected, and hence a topological ball. Removing a small $(B^5, B^3)$ pair with a standard action along the fixed-point set yields the desired concordance.

\end{proof}

\section{Examples of rank two group actions}

The \emph{standard linear} $G$ action is defined to be the restriction to the unit sphere in $\br^5$ of
$$ g\mapsto \begin{pmatrix}
\cos(2\pi/p)& -\sin(2\pi/p) & 0 & 0 & 0\\
\sin(2\pi/p) & \cos(2\pi/p) & 0 & 0 & 0\\
0& 0 & 1& 0& 0\\
0&0&0&1&0\\
0&0&0&0&1
\end{pmatrix}, h\mapsto \begin{pmatrix}
1&0&0&0&0\\
0&1&0&0&0\\
0&0&\cos(2\pi/p)& -\sin(2\pi/p) & 0\\
0&0&\sin(2\pi/p) & \cos(2\pi/p) & 0\\
0&0&0&0&1
\end{pmatrix}.$$

An easy exercise in representation theory shows that any orthogonal action (preserving orientation if $p=2$) is equivalent to the standard action modulo precomposition with an automorphism of $G$.

\subsection{A construction using counterexamples to the generalized Smith conjecture}

Giffen~\cite{Giffen} produced the first known examples of $\bz_p$ actions on $S^4$ with knotted two-sphere as fixed point set. He showed that, when $k\ge 2$ and $p=\pm 1\pmod{k}$, then the $p$-fold branched cover of any $k$-twist spun classical knot is a homotopy four-sphere, and diffeomorphic to $S^4$ when $p$ is odd. 

Gordon~\cite{Gordon2} re-cast Giffen's construction without the branched covering using a sort of equivariant Dehn surgery on the twist-spun knot. Gordon also gave another construction: If $(S^4, K)$ is any knot in $S^4$, and $g_p:S^4\to S^4$ is a standard rotation action on $S^4$ fixing $S^2$, then the equivariant connected sum of $p$ copies of $K$ admits a $\bz_p$ action (with unknotted fixed set). Removing a neighborhood of $pK$ and re-gluing it differently results in an action with knotted fixed set. 

Finally, Sumners\cite{Sumners} gave what is perhaps the simplest construction: Start with a standard rotation action on $B^5$ fixing $B^3$, glue on orbits of $1$-handles, and then glue in orbits of $2$-handles which cancel the $1$-handles geometrically in $\partial B^5$, but not even homotopically in $\partial(B^5\setminus B^3)$. The resulting action  on $\partial B^5$ has the nice additional property that a concordance to a standard linear action is obvious (and obviously smooth).

Actions of $\bz_p\times\bz_p$ on $S^4$ can be combined via equivariant connected sum at a fixed point. An action of $\bz_p\times\bz_p = \lan g, h\ran$ can also be modified by $\bz_p\times\bz_p$-equivariantly summing a $\lan g\ran$-orbit of $p$ copies of an $\lan h \ran$-action. The knot group of $\fix(h)$ is then an amalgamated product of nontrivial knot groups, so the action is non-linear.

\section{General analysis}

	We prove Theorem~\ref{maintheorem} by attempting to generalize the argument of Theorem~\ref{rankonecase} to construct a concordance, and measuring the obstructions. Henceforth, let $G=\bz_p\times\bz_p$.

In general, given a locally linear, homologically trivial  $G$ action on $S^4$, Smith Theory arguments show that $\fix(G)$ will consist of two points, 
$x_1$ and $x_2$. As the $G$-action is linear on $T_{x_1}S^4$, we can choose generators $g$ and $h$ of $G$ which fix transverse $2$-spheres, and so that after the two-spheres are oriented, $g$ and $h$ have
rotation numbers $(0,1)$ and $(1,0)$ at $x_1$. 
We label the spheres fixed by $g$ and $h$ $S_g$ and $S_h$, respectively.
Of course, $S_g\cap S_h =\{x_1, x_2\}$, while $S_g\cup S_h$, the
 entire singular
set of the action, will be denoted $\Sigma$.

We seek to analyze a given action via its quotient space, or more precisely, the
triple $(X, Y_g, Y_h) = (S^4/G, S_g/G, S_h/G)$. As in the case of Theorem~\ref{rankonecase}, $X$ is a simply-connected homology four-sphere, hence homeomorphic to $S^4$. 
Each of $Y_g$
and $Y_h$ is itself a two-sphere, over which $S_g$ and $S_h$ are $p$-fold cyclic
branched covers.  The original action can then be recovered
from the quotient as the $\bz_p\times\bz_p$  branched cover\footnote{more precisely, as an iterated branched cover, first by $\lan g\ran$, and then by $\lan h\ran$, or vice versa.} of $X$ over the knotted spheres $Y_g$ 
and $Y_h$. In this way, classification of  $G$-actions can be viewed in terms of knot 
theory.

As a model, let us denote   the quotient of $S^4$ by the standard action as $X_{\text{std}}$, and
choose  closed tubular neighborhoods $N(Y_g)$ and $N(Y_h)$ of $Y_g$ and $Y_h$ in $X_{\text{std}}$. Each is a trivial $D^2$-bundle over (an orbifold) $S^2$, and their union is a neighborhood $N(\overline{\Sigma})$ of the singular set $\overline{\Sigma}$ formed by 
plumbing the bundles together at $\overline{x}_1$ and $\overline{x}_2$. 
The closure of $X_{\text{std}}\setminus N(\Sigma)$ is homeomorphic to $D^2\times T^2$, with boundary $S^1\times T^2\homeo T^3$. The $D^2$ factor is bounded by a curve $\gamma_{\text{std}}$ running from $x_1$ to $x_2$ along a \lq\lq longitude" line of (a thickened) $Y_g$, then returning to $x_1$ along $Y_h$; finally, $D^2\times T^2$ has \lq\lq corners" along $\{x_1, x_2\}\times T^2$.

Our goal is to construct a cobordism rel boundary $W$ from $X\setminus N(\overline{\Sigma})$, the closure of the complement of a neighborhood of the singular set in $X=S^4/G$, to $X_{\text{std}}\setminus N(\overline{\Sigma})$, and then modify it so that $W\cup N(\overline{\Sigma})\times I\homeo S^4\times I$, and so that its iterated branched cover is also a cylinder. 

It follows from ~\cite[9.3]{FQ} that each of $Y_g$ and $Y_h$
has a topological normal bundle, so as before we can form $N(\overline{\Sigma})$. Alexander duality shows that $X\setminus N(\overline{\Sigma})$ has the homology of $D^2\times T^2$. Let $\mu_g$ and $\mu_h$ be small meridional loops in $\partial (X\setminus N(\overline{\Sigma}))$ around $Y_g$ and $Y_h$, respectively, and let $\gamma$ be a path running from $x_1$ to $x_2$ along $Y_g$, then back to $x_1$ along $Y_h$. Fix a homeomorphism $f$ from $\partial (X\setminus N(\overline{\Sigma}))$ to $\partial (X_{\text{std}}\setminus N(\overline{\Sigma}))$ sending $\mu_g, \mu_h,$ and $\gamma$ to their standard counterparts.

The space $D^2\times T^2$ is a $K(\bz\times\bz, 1)$, so the only potential obstruction to extending $f$ over all of $X\setminus N(\overline{\Sigma})$ lies in $H^2(X\setminus N(\overline{\Sigma}), \partial (X\setminus N(\overline{\Sigma}))); \bz\times\bz)$. This group is generated by the classes $\delta[\gamma_i^*]$, where $\gamma_1^*$ is a cocycle evaluating to $(1,0)$ on $\gamma$, and $\gamma_2^*$ is a cocycle evaluating to $(0,1)$.  As $f\circ \gamma = \gamma_{\text{std}}$ bounds a disk in $D^2\times T^2$, the obstruction to extending $f$ vanishes. Since $f$ is a homeomorphism on $\partial (X\setminus N(\overline{\Sigma}))$,  it is a degree one map relative to $\partial (X\setminus N(\overline{\Sigma}))$.  Finally, since $X\homeo S^4$, the stablized tangent bundle $\tau$ of $X\setminus N(\overline{\Sigma}))$ admits an obvious trivialization. Since $f^*(\varepsilon)$ and $\tau$ are individually trivial, 
trivializations can be combined to yield a stable trivialization $F$ of $\tau\oplus f^*\varepsilon$.
The data $((X\setminus N(\overline{\Sigma}), \partial), f, \varepsilon, F)$ together define a degree one normal map.

According to Freedman~\cite{FQ}, surgery \lq\lq works" for  four-manifolds with fundamental group $\bz\times\bz$, so there is an exact sequence\footnote{Since $Wh(\bz\times\bz)=0$, the distinction between homotopy and simple homotopy is not relevant to the structure set or
the  the surgery groups.}
$$ S_{\text{TOP}}(D^2\times T^2, \partial)\to NM_{\text{TOP}}(D^2\times T^2, \partial)\stackrel{\theta}{\to}L_4(\bz\times\bz).$$

$ S_{\text{TOP}}(D^2\times T^2, \partial)$ is a structure set measuring the difference between homotopy $D^2\times T^2$ manifolds with boundary $T^3$, and topological ones. The extension of the classification of \lq\lq fake tori" ~\cite[11.5]{FQ}  to dimension four shows that $ |S_{\text{TOP}}(D^2\times T^2, \partial)|=1$.

So $f$, together with the pullback $f^*\varepsilon$ of the trivial bundle, defines a surgery obstruction which vanishes if and only if $f$ is normally cobordant rel boundary to a homeomorphism. 

According to Wall ~\cite[13A.8, 13B.8]{Wall}, $L_4(\bz\times\bz)\cong\bz\oplus\bz_2$,
with generators for the factors given by the signature and a codimension two Kervaire-Arf invariant, respectively:  $\theta(f)=(\sigma(f), c(f))$. Here, $\sigma$ is the difference in signatures between $(X\setminus N(\overline{\Sigma}), \partial)$ and $(T^2\times D^2, \partial)$. Alexander duality shows that the two spaces have isomorphic cohomology, so the signature component of the obstruction vanishes. 

The second component, $c$ is more interesting. Composing the projection $p: D^2\times T^2\to T^2\homeo K(\bz\times\bz, 1)$ with $f$  realizes the obstruction as the 
element of $L_2(1)$ defined by the stable normal bundle of the transverse preimage of a regular value of $pf$. The submanifold $(pf)^{-1}(*)$ of $(X\setminus N(\overline{\Sigma}))$ is exactly the preimage via $f$ of a meridional disk, the first factor of $D^2\times T^2$. And then $c(f)$ is the Kervaire invariant of the associated normal map. 

%


%



This potential obstruction to null-concordance of our group action has the following equivalent interpretations: 
\begin{enumerate}
\item{The surgery obstruction (rel boundary), of the normal map from $X\setminus N(\Sigma)\to D^2\times T^2$.}
\item{The codimension two surgery obstruction of the restricted normal map $f^{-1}(D^2)\to D^2$}
\item{The element of the spin bordism group $\Omega_2^{\text{spin}}\cong\bz_2$ defined by identifying $-D^2$ to  $f^{-1}(D^2)$ and their normal bundles along their common boundary, using $F$.}
\end{enumerate}

We return to this invariant later. Let us suppose now that it vanishes. Then the map $f$ is normally cobordant to a homeomorphism. Let $(\overline{W}_0, f, F)$ be a normal cobordism. 

The cobordism $f: \overline{W}_0\to D^2\times T^2$ easily extends to a normal map $\overline{W}_0\to D^2\times T^2\times I$, which we continue to denote by $f$. 
Let $\mathcal{F}:\bz[\bz\times\bz]\to \bz[\bz_p\times\bz_p]$ be the homomorphism of group rings corresponding to the covering $S^4\setminus N(\Sigma)\to X\setminus N(\Sigma)$. The homology surgery groups of 
Cappell and Shaneson~\cite{CS2} measure the obstruction to surgering $f$ to be a $\bz[\bz_p\times\bz_p]$-homology equivalence. The relevant obstruction group $\Gamma_5^h(\mathcal{F})$ injects into $L_5(\bz_p\times\bz_p)$. The latter group vanishes for all primes $p$ (~\cite{WallSomeLGroups} is a convenient reference), so $\overline{W}_0$ is normally bordant rel boundary to a new bordism $\overline{W}_1$ which is $\bz[\bz_p\times\bz_p]$-homology equivalent to a cylinder. According to ~\cite[2.1]{CS2}, $f$ may be assumed to be 2-connected. The $\bz_p\times\bz_p$-cover $W_1$ of $\overline{W}_1$ therefore has the following properties:

\begin{enumerate}
\item{$\tilde{f}:W_1\to (S^4-N(\Sigma_{\text{std}}))\times I$ induces an isomorphism on integral homology.}
\item{$\partial{W_1}=S^4-N(\Sigma)\cup T^3\times I\cup S^4-N(\Sigma_{\text{std}})$.}
\item{$\pi_1(W_1)\cong \bz\times\bz$.}
\item{$W_1$ is equipped with a $\bz_p\times\bz_p$ action.}
\end{enumerate}

Equivariantly gluing in a copy of $N(\Sigma)\times I$ results in a simply-connected five-manifold homotopy equivalent, and hence homeomorphic, to $S^4\times I$, and hence yields the desired concordance.

\section{Behavior of the obstruction}

Are there examples of $\bz_p\times\bz_p$-actions which are not concordant to linear actions? We do not yet have a definitive answer to this question, but as we shall see, the obstruction $c$ vanishes for known constructions of $\bz_p\times\bz_p$-actions. 

\begin{eg}
Montesinos~\cite{MontesinosTwins1, MontesinosTwins2} considered \lq\lq twins": pairs of two-spheres in $S^4$ which intersect transversely in two points. Any such pair of twins has a regular neighborhood homeomorphic to a pair of trivial $D^2$ bundles over $S^2$ plumbed together at their north and south poles. As before, the obstruction to extending a homeomorphism on the boundaries of these regular neighborhoods to a map of the twin complements to $D^2\times T^2$ vanishes, so our definition of the Arf invariant in fact extends to the complement of any pair of \lq\lq twin" two-spheres.

Pairs of twins arise naturally in the spinning of classical knots: One sphere is formed by the spun knot, and the other is the boundary of the ball $B^3$, held fixed during the spinning. A map  $f:B^3\setminus B^1\to S^1$ defining a Seifert surface relative to the boundary of the spinning ball extends to 
$$f\times{\text{id}}:(B^3\setminus B^1)\times S^1\to S^1\times S^1$$
so that the same Seifert surface used to calculate the Arf invariant of the classical knot can also be used to calculate $c$ (cf.~\cite{Robertello}). So, for example, the pair of twins (Spun trefoil, $\partial B^3$) realizes $c=1\in (\bz_2, +)$.  Note that the branched covers of $S^4$ over the spun trefoil are $S^1\times\bz_p$-manifolds, but by the solution of the classical Smith conjecture, they are never four-spheres. 

\end{eg}

\begin{prop} The invariant $c$ is additive under connected sum of $\bz_p\times\bz_p$ actions at a fixed point, and is unchanged when a given $\bz_p\times\bz_p$ action is modified by equivariant connected sum with $p$ copies of a $\bz_p$ action.
\end{prop}

\begin{proof}

First consider the case of two actions joined at a fixed point. Note that $c$ is calculated by extending a standard map $f:\partial N(\Sigma)\to T^2$ to $X\setminus N(\Sigma)$, then examining the Kervaire-Arf invariant on the normal bundle of  $f^{-1}(*)$. Because the map and the action are standard on  $\partial N(\Sigma)$, the maps ($f_1$ and $f_2$, say) may be assumed to agree where their domains overlap, so that the new surface used to compute $c$ is the boundary connected sum $f_1^{-1}(*)\natural f_2^{-1}(*)$. As in the case of Seifert surfaces for classical knots, the Arf invariant is additive under this operation.

Now consider the case of a $G$-action summed with a $\lan g\ran$-orbit of $\lan h\ran$-actions. 
We claim that an appropriate map $f:S^4\setminus \Sigma\to T^2$ may be modified to yield a new map for the revised action by changing $f$ only near the sum points. Since we may simply choose a point $y\in T^2$ away from the image of these sum points to calculate $c$, it will follow immediately that $c$ is unchanged.

As before, the actions may be assumed standard on an equivariant regular neighborhood of the singular set.  A map $f:S^4\setminus N(\Sigma)\to T^2$ which is convenient for the calculation of $c$ is obtained as follows: as the quotient $(S^4/\lan g \ran, Y_g)$ is again a knotted sphere pair, the complement $S^4/\lan g\ran\setminus Y_g$ admits a map $\overline{f}_g$ to $S^1$ inducing an isomorphism on $H_1$. Then $\overline{f}_g$ lifts to an equivariant map $f_g:S^4\setminus S_g\to\widetilde{S}^1\cong S^1$. A similar construction applies to $S_4\setminus S_h$. We define $f=(f_g\times f_h)|_{S^4\setminus\Sigma}:S^4\setminus\Sigma \to S^1\times S^1$.

Let $x$ be a point on $S_h\subset S^4$, away from $\fix(g)$, such that $f(x)\in T^2$ is a regular value for $f$.  As the $G$-action is standard near $\Sigma$,  a neighborhood of $x$ may be parametrized
in the form $B^2\times (-\delta, \delta)\times (-\epsilon, \epsilon)$, so that $h\cdot(re^{i\theta}, t_1, t_2)=(re^{i(\theta+\frac{2\pi}{p})}, t_1, t_2)$, and so that, near $f(x)\in T^2$, $f(re^{i\theta}, t_1, t_2)=(e^{i\theta},t_2)$.

Meanwhile, let $S^{4'}$ be a sphere with an $\lan h\ran$-action fixing $S_h'$, and let $f_h':S^4;\setminus S_h\to S^1$ be an equivariant Seifert map. Choose a point $x'\in S_h'$; near $x'$, choose coordinates so that 
$f_h'(re^{i\theta}, t_1, t_2) =e^{i\theta}$, and extend $f_h' $ to $S^1\times (-\epsilon, \epsilon)$ in the obvious way, so that $p$ copies of $S^4;$ may be glued $G$-equivariantly to $\lan g\ran\times \{ x\}$ in $S^4$. 

\end{proof}

\section{Concluding remarks and questions}

It follows from the additivity of the Arf invariant that the constructions of $\bz_p\times\bz_p$ actions outlined in the previous section yield actions which are concordant to linear actions. Modifications to these actions to change the Arf invariant would need to be sufficiently global as to change $f^{-1}(x)$ for every $x\in T^2$.


It is not clear whether such \lq\lq highly nonlinear" actions exist, but it is interesting to compare the classical case of $\bz_2$ actions in dimension three. Murasugi~\cite{MurasugiArf} proved that
for a knot $K$, $\arf(K) =0 \iff \Delta_K(-1)\cong \pm 1 \pmod 8$. As $|\Delta_K(-1)|$ is the order of the first homology of the two-fold  branched cover of $K$, $\arf(K)$ is therefore an obstruction to a knot having a homology sphere as a double branched cover, and hence also an obstruction (which does not require the full strength of the classical Smith Conjecture) to the existence of a $\bz_2$-action fixing $K$. 
A simplified, essentially homological proof of this fact might  generalize to yield an obstruction for  $\bz_2\times\bz_2$ actions on $S^4$. The realization question for odd primes seems harder. We hope to consider both in future work.

\bibliographystyle{abbrv}
\bibliography{mybiblio}

 \end{document}